\documentclass[12pt]{amsart}

\usepackage{fullpage}
\usepackage{amsmath}
\usepackage{amssymb}
\usepackage{theorem}
\usepackage{hyperref}
\usepackage{units}
\usepackage{graphicx, pstcol, pst-plot, pst-node}

\newcommand{\sect}{\section}
\newcommand{\bi}{\begin{itemize}}
\newcommand{\ei}{\end{itemize}}

\newcommand{\Av}{\operatorname{Av}}
\newcommand{\gr}{\operatorname{gr}}
\newcommand{\conv}{\operatorname{conv}}
\newcommand{\ugr}{\operatorname{\overline{gr}}}
\newcommand{\lgr}{\operatorname{\underline{gr}}}
\newcommand{\proble}{\bigskip\noindent{\bf Problem:} }
\newcommand{\problem}[1]{\bigskip\noindent{\bf Problem (#1):} }
\newcommand{\problemm}[2]{\bigskip\noindent{\bf \$#1 Problem (#2):} }
\newcommand{\conjecture}[1]{\bigskip\noindent{\bf Conjecture (#1):} }
\newcommand{\conjecturem}[2]{\bigskip\noindent{\bf \$#1 Conjecture (#2):} }
\newcommand{\background}{\bigskip\noindent{\bf Background:} }

\begin{document}

\title{Problems and Conjectures presented at the Third International Conference on Permutation Patterns\\ (University of Florida, March 7--11, 2005)}
\author{Recounted by Murray Elder and Vince Vatter}

\maketitle

\sect{A very brief introduction to permutation patterns}

A {\em permutation} is an arrangement of a finite number of distinct elements of a linear order, for example, $e,\pi,0,\sqrt 2$ and $3412$. Two permutations are {\em order isomorphic} if the have the same relative ordering.
We say a permutation $\tau$ {\em contains} or {\em involves} a permutation $\beta$ if deleting some of the entries of $\pi$ gives a permutation that is order isomorphic to $\beta$, and we write $\beta\le\tau$. For example, $534162$ (when permutations contain only single digit natural numbers we suppress the commas) contains $321$ (delete the values $4$, $6$, and $2$). A permutation {\em avoids} a permutation if it does not contain it.

For a set of permutations $B$ define $\Av(B)$ to be the set of permutations that avoid all of the permutations in $B$ and let $s_n(B)$ denote the number of permutations of length $n$ in $\Av(B)$.  A set of permutations or class $C$ is {\em closed} if $\pi\in C$ and $\sigma\leq \pi$ implies $\sigma\in C$.  Therefore $\Av(B)$ is a class for every set of permutations $B$ and every permutation class can be written as $\Av(B)$ for some set $B$.   An {\em antichain} of permutations is a set of permutations such that no permutation contains another.
We call a set of permutations $B$ the {\em basis} of a class of permutations $C$ if $C=\Av(B)$ and $B$ is an antichain.

\bigskip

\noindent{\bf For more information:}
\bi
\item  B\'ona's book~\cite{bona:book} gives a good introduction to permutation patterns, 
\item Kitaev and Mansour have a survey~\cite{km:survey} on the topic,
\item Volume 9 (2) of the {\it Electronic Journal of Combinatorics\/} is devoted to permutation patterns,
\item {\it Advances in Applied Mathematics\/} will have a special issue on permutation patterns soon.
\ei

\sect{Exact enumeration}

\background A sequence $s_n$ is said to be $P$-recursive (which is short for {\em polynomially recursive\/}) if there 
are polynomials $p_0,p_1,\ldots p_k$ so that 
$$p_k(n)s_{n+k}+p_{k-1}(n)s_{n+k-1}+\ldots +p_0(n)s_n=0$$
For example, $s_n=n!$ is $P$-recursive since $s_{n+1}-(n+1)s_n=0$.  A sequence is $P$-recursive if and only if its generating function is $D$-finite, meaning that its derivatives span a finite dimensional vector space over $\mathbb{C}(x)$.  More information on $P$-recursive sequences can be found in Stanley~\cite{stanley:prec} and Zeilberger~\cite{z:holo1}.

Let $\pi_1,\pi_2,\dots,\pi_k$ be permutations and $r_1,r_2,\dots,r_k$ positive integers.  Noonan and Zeilberger~\cite{nz:copies} conjectured that the number of $n$-permutations with exactly $r_i$ copies of $\pi_i$ for each $i\in[k]$ forms a $P$-recursive sequence in $n$.  Atkinson~\cite{a:rp} showed that this is equivalent to the seemingly weaker conjecture\footnote{Some have called this Gessel's Conjecture, although to be fair to Gessel, all he wrote in \cite{gessel} was that \begin{quote}Another possible candidate for $P$-recursiveness is the problem of counting permutations (or more generally sequences) with forbidden subsequences defined by inequalities, for example permutations $a_1a_2\cdots a_n$ of $\{1,2,\dots,n\}$ with no subsequence $a_ia_ja_k$ satisfying $a_i<a_k<a_j$.\end{quote}} that $s_n(B)$ is $P$-recursive for all finite sets of patterns $B$.

For $B=\{1324\}$, Marinov and Radoi{\v{c}}i{\'c}~\cite{mr:1324} used a generating tree approach to find the first 20 terms.
Recent work by Albert, Elder, Rechnitzer, Westcott and Zabrocki~\cite{1324} gives more terms for the sequence $s_n(1324)$ (another six or so) which suggest this sequence is not well behaved.

\conjecture{Doron Zeilberger} The Noonan-Zeilberger Conjecture is false, and, in particular, the sequence $s_n(1324)$ is not $P$-recursive.
\bigskip

In fact, Zeilberger made the stronger claim that ``not even God knows $a_{1000}(1324)$.''  He\footnote{Zeilberger, not God} suggests that we replace the Noonan-Zeilberger Conjecture with the following research program:

\problem{Doron Zeilberger} Find necessary and sufficient conditions characterizing the classes whose generating functions are:
\bi
\item rational,
\item algebraic,
\item $D$-finite,
\item $\dots$
\ei
Currently we have only one result of this type: Kaiser and Klazar~\cite{kk:growth} characterize the classes with polynomial enumeration.

One can also consider pattern avoiding compositions of an integer.  Let $${\mathbf a}=(a_1,a_2,\dots,a_k)$$ denote a composition of $n$, so $a_1+a_2+\cdots+a_k=n$, the $a_i$s need not be distinct, and order matters.  We say that ${\mathbf a}$ contains the permutation $\beta$ if it contains a subsequence that is order isomorphic to $\beta$.

\conjecturem{201}{Herb Wilf} The number of $123$-avoiding compositions of $n$ into positive parts is not a $P$-recursive sequence.

Savage and Wilf~\cite{sw:comp} have found the generating function for this sequence:
$$
\sum_{i\ge 1}\frac{1}{1-x^i}\prod_{j\neq i}\left\{\frac{1-x^i}{(1-x^{j-i})(1-x^i-x^j)}\right\}.
$$
Zeilberger pointed out that this function appears to have infinitely many singularities, indicating that it is not $D$-finite.

\sect{Comparison of patterns}

\background We say that the sets $B_1$ and $B_2$ are {\it Wilf-equivalent\/}, and write $B_1\sim B_2$, if $s_n(B_1)=s_n(B_1)$ for all natural numbers $n$, that is, if $B_1$ and $B_2$ are equally avoided.  For example, it is a classic result that every permutation in $S_3$ is Wilf-equivalent to every other permutation in $S_3$, or in other words, $S_3$ contains only one Wilf-equivalence class.  It is a much more recent result (due to Stankova~\cite{stankova:len4, stankova:fs} and West~\cite{west:cat}) that $S_4$ contains three Wilf-equivalence classes.  Although it was not raised during the problem session, we can not resist mentioning the following tantalizing problem.

\proble Find necessary and sufficient conditions for two permutations to be Wilf-equivalent.

\bigskip

Many had thought that the permutations of a certain length could be ordered according to avoidance, modulo Wilf-equivalence.  Stankova and West~\cite{sw:wilf} were the first to find a counterexample.  They observed that $s_n(53241)<s_n(43251)$ for all $n\le 12$ while $s_{13}(53241)>s_{13}(43251)$.  Still, they conjectured that the permutations of a certain length could be asymptotically ordered according to avoidance, modulo Wilf-equivalence, where we say that $\pi$ is asymptotically more avoidable that $\sigma$ if $s_k(\pi)>s_k(\sigma)$ for all sufficiently large $k$.

The conjectures in this section both attempt to understand these orderings, but before stating them we need some additional notation.

\begin{figure}
\psset{xunit=0.00714286in, yunit=0.00714286in}
\psset{linewidth=0.6\psxunit}
\begin{pspicture}(0,0)(150,150)
\psframe[linewidth=0.01in](0,0)(150,150)
\psline[linecolor=gray, linewidth=2.0\psxunit](5,5)(5,85)
\psline[linecolor=gray, linewidth=2.0\psxunit](5,85)(85,85)
\psline[linecolor=gray, linewidth=2.0\psxunit](85,85)(85,5)
\psline[linecolor=gray, linewidth=2.0\psxunit](85,5)(5,5)
\psline[linecolor=gray, linewidth=2.0\psxunit](85,85)(85,145)
\psline[linecolor=gray, linewidth=2.0\psxunit](85,145)(145,145)
\psline[linecolor=gray, linewidth=2.0\psxunit](145,145)(145,85)
\psline[linecolor=gray, linewidth=2.0\psxunit](145,85)(85,85)
\pscircle*(10,20){3\psxunit}
\pscircle*(20,40){3\psxunit}
\pscircle*(30,60){3\psxunit}
\pscircle*(40,80){3\psxunit}
\pscircle*(50,10){3\psxunit}
\pscircle*(60,30){3\psxunit}
\pscircle*(70,50){3\psxunit}
\pscircle*(80,70){3\psxunit}
\pscircle*(90,110){3\psxunit}
\pscircle*(100,90){3\psxunit}
\pscircle*(110,130){3\psxunit}
\pscircle*(120,100){3\psxunit}
\pscircle*(130,140){3\psxunit}
\pscircle*(140,120){3\psxunit}
\end{pspicture}
\caption{The plot of $24681357\oplus 315264=2,4,6,8,1,3,5,7,11,9,13,10,14,12$.}\label{exsum}
\end{figure}
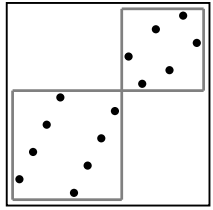

Given two permutations $\pi\in S_m$ and $\sigma\in S_n$, we define the {\em direct sum} (or simply {\em sum}) of $\pi$ and $\sigma$, written $\pi\oplus\sigma$ by
$$
\pi\oplus\sigma(i)
=
\left\{
\begin{array}{ll}
\pi(i)&\mbox{if $i\in [m]$,}\\
\sigma(i-m)&\mbox{if $i\in[m+n]\setminus[m]$.}
\end{array}
\right.
$$
Figure~\ref{exsum} shows an example.

A permutation is said to be {\em layered} if it is the direct sum of some number of decreasing permutations.  For example, $21\oplus 321\oplus 1=215436$ is a layered permutation.

\begin{figure}
\begin{center}
\psset{xunit=0.00714286in, yunit=0.00714286in}
\psset{linewidth=0.6\psxunit}
\begin{pspicture}(0,0)(150,150)
\psframe[linewidth=0.01in](0,0)(150,150)
\psline[linecolor=gray, linewidth=2.0\psxunit](5,95)(5,115)
\psline[linecolor=gray, linewidth=2.0\psxunit](5,115)(25,115)
\psline[linecolor=gray, linewidth=2.0\psxunit](25,115)(25,95)
\psline[linecolor=gray, linewidth=2.0\psxunit](25,95)(5,95)
\psline[linecolor=gray, linewidth=2.0\psxunit](25,5)(25,65)
\psline[linecolor=gray, linewidth=2.0\psxunit](25,65)(85,65)
\psline[linecolor=gray, linewidth=2.0\psxunit](85,65)(85,5)
\psline[linecolor=gray, linewidth=2.0\psxunit](85,5)(25,5)
\psline[linecolor=gray, linewidth=2.0\psxunit](85,115)(85,145)
\psline[linecolor=gray, linewidth=2.0\psxunit](85,145)(115,145)
\psline[linecolor=gray, linewidth=2.0\psxunit](115,145)(115,115)
\psline[linecolor=gray, linewidth=2.0\psxunit](115,115)(85,115)
\psline[linecolor=gray, linewidth=2.0\psxunit](115,65)(115,95)
\psline[linecolor=gray, linewidth=2.0\psxunit](115,95)(145,95)
\psline[linecolor=gray, linewidth=2.0\psxunit](145,95)(145,65)
\psline[linecolor=gray, linewidth=2.0\psxunit](145,65)(115,65)
\pscircle*(10,100){3\psxunit}
\pscircle*(20,110){3\psxunit}
\pscircle*(30,30){3\psxunit}
\pscircle*(40,10){3\psxunit}
\pscircle*(50,50){3\psxunit}
\pscircle*(60,20){3\psxunit}
\pscircle*(70,60){3\psxunit}
\pscircle*(80,40){3\psxunit}
\pscircle*(90,130){3\psxunit}
\pscircle*(100,140){3\psxunit}
\pscircle*(110,120){3\psxunit}
\pscircle*(120,90){3\psxunit}
\pscircle*(130,80){3\psxunit}
\pscircle*(140,70){3\psxunit}
\end{pspicture}
\end{center}
\caption{The plot of $3142[12,315264,231,321]=10,11,3,1,5,2,6,4,13,14,12,9,8,7$.}\label{exwreath}
\end{figure}
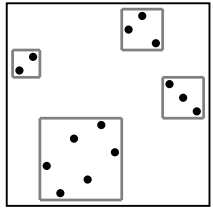

Let $\pi$ be a permutation of length $k$ and let $\sigma_1,\sigma_2,\dots,\sigma_k$ also be permutations.  We define the {\it inflation\/} of $\pi$ by $\sigma_1,\sigma_2,\dots,\sigma_k$ to be the permutation obtained by replacing the $i$th element of $\pi$ by a block order isomorphic to $\sigma_i$ so that the blocks themselves are order isomorphic to $\pi$.  We write the resulting permutation as $\pi[\sigma_1,\sigma_2,\dots,\sigma_k]$.  An example is shown in Figure~\ref{exwreath}.  Direct sums fit into this framework as well, since $\pi\oplus\sigma=12[\pi,\sigma]$.

\conjecture{Alex Burstein} For any permutations $\sigma_1, \sigma_2, \sigma_3$, and any natural number $n$,
\begin{enumerate}
\item $s_n(132[\sigma_1,\sigma_2, \sigma_3]) < s_n(123[\sigma_1,\sigma_2, \sigma_3])$,
\item $s_n(312[\sigma_1,\sigma_2, \sigma_3]) < s_n(123[\sigma_1,\sigma_2, \sigma_3])$,
\item $s_n(123[\sigma_1,12\dots t, \sigma_3]) < s_n(123[\sigma_1,t\dots 21, \sigma_3])$ for any $t\ge 2$.
\end{enumerate}

\conjecture{Alex Burstein} For any nonlayered permutation $\pi\in S_k$, permutations $\sigma_1,\sigma_2,\dots,\sigma_k$, and natural numbers $n$,
$$
s_n(\pi[\sigma_1,\sigma_2,\dots,\sigma_k]) <
s_n(12\dots k[\sigma_1,\sigma_2,\dots,\sigma_k]).
$$

Of course, we don't have a proof of this latter conjecture even in the special case $\sigma_1=\sigma_2=\dots=\sigma_k=1$, and if we did, we would have a different proof of the Stanley-Wilf Conjecture via Regev~\cite{regev}.

\begin{figure}[ht!]
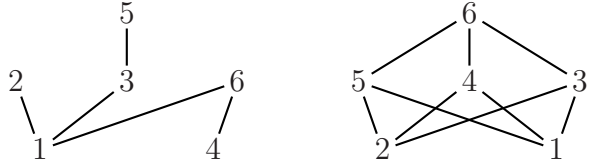

\begin{tabular}{ccc}
\psset{nodesep=2pt,colsep=18pt,rowsep=18pt}
\begin{psmatrix}
&&[name=5] 5&&\\[0pt]
[name=2]2&&[name=3]3&&[name=6]6\\[0pt]
[name=1]1\psspan{2}&&[name=4]4\psspan{2}
\ncline{5}{3}
\ncline{2}{1}
\ncline{3}{1}
\ncline{6}{1}
\ncline{6}{4}
\end{psmatrix}
&\rule{20pt}{0pt}&
\psset{nodesep=2pt,colsep=18pt,rowsep=18pt}
\begin{psmatrix}
&&[name=6]6&&\\[0pt]
[name=5]5&&[name=4]4&&[name=3]3\\[0pt]
[name=2]2\psspan{2}&&[name=1]1\psspan{2}
\ncline{6}{3}
\ncline{6}{4}
\ncline{6}{5}
\ncline{3}{1}
\ncline{3}{2}
\ncline{4}{1}
\ncline{4}{2}
\ncline{5}{1}
\ncline{5}{2}
\end{psmatrix}
\end{tabular}
\caption{The posets corresponding to $416352$ (left) and $215436$ (right)}\label{poset-first}
\end{figure}

B\'ona's conjecture takes a little more preparation.  To any permutation $\pi$ on $[n]$ we associate a poset $P_\pi$.  The elements of this poset are the integers $1$ through $n$ (which we will think of as the values of $\pi$, although we could equivalently think of them as being the indices of $\pi$).  We will denote this partial order by $\preceq$, and write $i\preceq j$ if $i\le j$ (as integers) and if $i$ occurs before $j$ in $\pi$ (reading from left to right).  Figure~\ref{poset-first} shows two examples.

\begin{figure}[ht!]
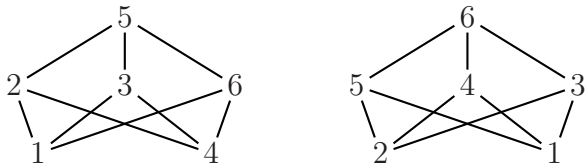

\begin{tabular}{ccc}
\psset{nodesep=2pt,colsep=18pt,rowsep=18pt}
\begin{psmatrix}
&&[name=5] 5&&\\[0pt]
[name=2]2&&[name=3]3&&[name=6]6\\[0pt]
[name=1]1\psspan{2}&&[name=4]4\psspan{2}
\ncline{5}{3}
\ncline{5}{2}
\ncline{5}{6}
\ncline{2}{4}
\ncline{3}{4}
\ncline{2}{1}
\ncline{3}{1}
\ncline{6}{1}
\ncline{6}{4}
\end{psmatrix}
&\rule{20pt}{0pt}&
\psset{nodesep=2pt,colsep=18pt,rowsep=18pt}
\begin{psmatrix}
&&[name=6]6&&\\[0pt]
[name=5]5&&[name=4]4&&[name=3]3\\[0pt]
[name=2]2\psspan{2}&&[name=1]1\psspan{2}
\ncline{6}{3}
\ncline{6}{4}
\ncline{6}{5}
\ncline{3}{1}
\ncline{3}{2}
\ncline{4}{1}
\ncline{4}{2}
\ncline{5}{1}
\ncline{5}{2}
\end{psmatrix}
\end{tabular}
\caption{The convex hulls of the posets from Figure~\ref{poset-first}}\label{poset-second}
\end{figure}

These posets are ranked, and the rank of $i$ in $P_{\pi}$ is the length of the longest
increasing subsequence of $\pi$ that ends in $i$.  Given a poset $P$
with rank function $r$, we let $\conv(P)$ denote the {\em convex hull}
of $P$.  The poset $\conv(P)$ is defined on the same set of elements
as $P$, but in $\conv(P)$, $i\preceq_{\conv} j$ if and only if
$r(i)\le r(j)$ (in $P$).  Figure~\ref{poset-second} shows the convex
hulls of the posets from Figure~\ref{poset-first}.  Note that
$\conv(P_{215436})=P_{215436}$ and that $\conv(P_{416352})\cong
P_{215436}$.  Indeed, these are not accidents.  For any layered
permutation $\pi$, $\conv(P_\pi)=P$, and for any permutation $\pi$
there is a unique layered permutation, which we refer to as $\conv(\pi)$,
such that $P_\pi\cong P_{\conv(\pi)}$.

\bigskip

\conjecture{Mikl\'os B\'ona} For any permutation $\pi$ and natural number $n$, $s_n(\pi)\le s_n(\conv(\pi))$.

Like Burstein's Conjecture, B\'ona's Conjecture would give another proof of the Stanley-Wilf Conjecture, via B\'ona~\cite{bona:layered}.

\sect{Growth rates}

\background Let $X$ be a set of permutations of varying lengths.  We define the {\it upper growth rate of $X$\/} to be
$$
\ugr(X)=\limsup_{n\rightarrow\infty}\sqrt[n]{|X\cap S_n|}
$$
and the {\it lower growth rate of $X$\/} to be
$$
\lgr(X)=\liminf_{n\rightarrow\infty}\sqrt[n]{|X\cap S_n|}.
$$

The Marcus-Tardos Theorem~\cite{mt:swc}, formerly the Stanley-Wilf Conjecture, implies that these limits are finite if $X$ is a closed class and $X$ is not the set of all permutations.  For this reason some refer to $\gr(\Av(B))$ as the ``Stanley-Wilf limit of $\Av(B)$.''

Of course it would be natural to define the {\it growth rate of $X$\/} as
$$
\gr(X)=\lim_{n\rightarrow\infty}\sqrt[n]{|X\cap S_n|},
$$
but this limit is not not known to exist in general.  In one important special case, Arratia~\cite{arratia} used Fekete's Lemma to show that $\gr(\Av(\pi))$ always exists, but his argument does not generalize to the case of $\Av(B)$ where $|B|\ge 2$.  For any $\pi\in S_3$, the growth rate of $\Av(\pi)$ is 4, because the $\pi$-avoiding permutations are counted by the Catalan numbers.  The growth rate of $\Av(12\dots k)$ is $(k-1)^2$ by Regev~\cite{regev}.  It was an old conjecture that the growth rate of $\Av(\pi)$ is $(k-1)^2$ for all $\pi\in S_k$, but this was disproved by B\'ona~\cite{b:1342}, who showed that the growth rate of $\Av(1342)$ is $8$.  It is also tempting from this data to conjecture that the growth rate of $\Av(\pi)$ is always an integer; B\'ona~\cite{bona:notint} disproved this also by showing that the growth rate of $\Av(12453)$ is $9+4\sqrt{2}$.

The only permutation of length four for which the growth rate of $\Av(\pi)$ is not known is $1324$.  Albert, Elder, Rechnitzer, Westcott, Zabrocki~\cite{1324} recently showed that this growth rate is at least 9.35, thereby disproving a conjecture of Arratia.\footnote{Arratia~\cite{arratia} had conjectured the growth rate of $\Av(\pi)$ is at most $(k-1)^2$ for any $\pi\in S_k$.}  AERWZ postulate the actual growth rate for $\Av(\pi)$ lies between 11 and 12.  This brings us to

\problemm{200}{Doron Zeilberger} Improve on the AERWZ lower bound for $\gr(\Av(1324))$ using the Maple package \texttt{WILFPLUS}~\cite{wilfplus:prog}.

\conjecture{Vince Vatter} For any finite set of patterns $B$, the growth rate of $\Av(B)$ exists.

\problem{Michael Albert} How quickly can an antichain of permutations grow?  More precisely, how large can $\ugr(X)$ be for an antichain $X$?

Albert has an example of an antichain for which this limit is approximately $2.24$.

\sect{Stack sorting}

\background A stack is a last-in first-out linear sorting device with push and pop operations.  The greedy algorithm for stack sorting a permutation $\pi=\pi(1)\pi(2)\dots\pi(n)$ goes as follows.  First we push $\pi(1)$ onto the stack.  Now suppose at some later stage that the letters $\pi(1),\dots,\pi(i-1)$ have all been either output or pushed on the stack, so we are reading $\pi(i)$.  We push $\pi(i)$ onto the stack if and only $\pi(i)$ is lesser than any element on the stack.  Otherwise we pop elements off the stack until $\pi(i)$ is less than any remaining stack element and we push $\pi(i)$ onto the stack.  This produces a permutation $s(\pi)$.  A permutation is {\it West $t$-stack sortable\/} if $s^t(\pi)$ is the identity permutation.  

Of course other stack sorting algorithms are possible, and for more than one stack, the iterated greedy algorithm described above is not the optimal algorithm (see, for example, Smith~\cite{smith:sorting}).  By allowing any sorting algorithm we reach the definition of {\it general $t$-stack sortability\/}.  In addition to the other general references on permutation patterns, B\'ona~\cite{bona:stacksurvey} provides a survey of stack sorting.

This more general sorting machine has ``complete lookahead,'' meaning that it can make its choices based on looking at the entire permutation.  Julian West suggested that it might be interesting to study machines that can only see the next $k$ letters of the permutations at a time.

\problem{Mikl\'os B\'ona} For which $t,n$ is the number of West $t$-stack sortable permutations of length $n$ an odd number?  So far we know
\bi
\item one stack: odd if and only if $n=2^k-1$ since 1-stack sortable permutations are counted by the Catalan numbers,
\item two stacks: odd quite frequently; more, precisely, this number is odd if and only if the binary expansion of $n$ does not contain two $1$'s in consecutive positions and ends in a $1$; the number of such values of $n$ between $1$ and $2^m$ is the $m$th Fibonacci number.  B\'ona~\cite{bona:parity} gives a combinatorial proof,
\item three stacks: numerical eveidence of Oliver Guibert: odd for n=1,9.
\ei
Thus with one stack odd numbers show up very rarely, with two stacks the numbers are quite often odd, but for three stacks odd numbers are again rare.  B\'ona conjectures that this pattern continues.  He provided the following motivation: $1$-stack sortable $n$-permutations can be described by a codeword on two letters of length $2n$ (using one letter for ``in'' and another for ``out''); $2$-stack sortable $n$-permutations can be described by a codeword on three letters of length $3n$ (using letters for ``into first stack,'' ``from first to second stack,'' and ``out of second stack'').  In the same manner, $t$-stack sortable $n$-permutations can be described by a codeword on $t+1$ letters of length $(t+1)n$.

Bruce Sagan pointed out that there may be a connection to trees that would help establish B\'ona's conjecture.

\problem{Mik\'os B\'ona} In particular, the number of West $2$-stack sortable permutations on $4k+3$ elements is always even.  Give a combinatorial proof.

\problem{Mike Atkinson} Can one decide in polynomial time whether an input permutation is sortable by a general $t$-stack sorting machine, or is this problem NP-complete?

Note that this problem does lie in NP, because one can verify in polynomial time whether a permutation can be sorted in some particular set of steps.  If the basis for general $2$-stack sortable permutations were finite, then Atkinson's problem could be done in polynomial time, but Murphy showed in his thesis~\cite{maximillian} that this basis is infinite.

\problem{Steve Waton} What is the shortest permutation that can not be sorted by a general 3-stack sorting machine?

Tarjan~\cite{tarjan:sorting} gave upper and lower bounds for the shortest permutation that can't be sorted by a $t$-stack machine, but they don't appear to be sharp.  Murphy~\cite{maximillian} shows that the shortest permutations that can't be be sorted by a general $2$-stack sorting procedure have length $7$, and there are $22$ of them.  For three stacks, Elder and Waton have wagered a beer on the problem, with Elder guessing that the shortest unsortable permutation is of length $15$  and Waton betting on $22$.

\bigskip
\bibliographystyle{acm}
\begin{small}
\bibliography{probrefs}

\begin{thebibliography}{10}

\bibitem{1324}
{\sc Albert, M., Elder, M., Rechnitzer, A., Westcott, P., and Zabrocki, M.}
\newblock On the {W}ilf-{S}tanley limit of $4231$-avoiding permutations and a
  conjecture of arratia.
\newblock
  \href{http://front.math.ucdavis.edu/math.CO/0502504}{arXiv:math.CO/0502504}.

\bibitem{arratia}
{\sc Arratia, R.}
\newblock On the {S}tanley-{W}ilf conjecture for the number of permutations
  avoiding a given pattern.
\newblock {\em Electron. J. Combin. 6\/} (1999), Note, N1, 4 pp. (electronic).

\bibitem{a:rp}
{\sc Atkinson, M.~D.}
\newblock Restricted permutations.
\newblock {\em Discrete Math. 195}, 1-3 (1999), 27--38.

\bibitem{bona:notint}
{\sc B{\'o}na, M.}
\newblock The limit of a {S}tanley-{W}ilf sequence is not always an integer!
\newblock {\it J. Comb. Theory, Ser. A\/}, to appear.

\bibitem{bona:parity}
{\sc B\'ona, M.}
\newblock Parity and stack sortability.
\newblock in preparation.

\bibitem{b:1342}
{\sc B{\'o}na, M.}
\newblock Exact enumeration of {$1342$}-avoiding permutations: a close link
  with labeled trees and planar maps.
\newblock {\em J. Combin. Theory Ser. A 80}, 2 (1997), 257--272.

\bibitem{bona:layered}
{\sc B{\'o}na, M.}
\newblock The solution of a conjecture of {S}tanley and {W}ilf for all layered
  patterns.
\newblock {\em J. Combin. Theory Ser. A 85}, 1 (1999), 96--104.

\bibitem{bona:stacksurvey}
{\sc B{\'o}na, M.}
\newblock A survey of stack-sorting disciplines.
\newblock {\em Electron. J. Combin. 9}, 2 (2002/03), Article 1, 16 pp.
  (electronic).

\bibitem{bona:book}
{\sc B{\'o}na, M.}
\newblock {\em Combinatorics of permutations}.
\newblock Discrete Mathematics and its Applications (Boca Raton). Chapman \&
  Hall/CRC, Boca Raton, FL, 2004.

\bibitem{gessel}
{\sc Gessel, I.~M.}
\newblock Symmetric functions and {$P$}-recursiveness.
\newblock {\em J. Combin. Theory Ser. A 53}, 2 (1990), 257--285.

\bibitem{kk:growth}
{\sc Kaiser, T., and Klazar, M.}
\newblock On growth rates of closed permutation classes.
\newblock {\em Electron. J. Combin. 9}, 2 (2002/03), Research paper 10, 20 pp.
  (electronic).

\bibitem{km:survey}
{\sc Kitaev, S., and Mansour, T.}
\newblock A survey on certain pattern problems.
\newblock University of Kentucky research report 2003--09, available at
  \url{http://www.ms.uky.edu/~math/MAreport/survey.ps}.

\bibitem{mt:swc}
{\sc Marcus, A., and Tardos, G.}
\newblock Excluded permutation matrices and the {S}tanley-{W}ilf conjecture.
\newblock {\em J. Combin. Theory Ser. A 107}, 1 (2004), 153--160.

\bibitem{mr:1324}
{\sc Marinov, D., and Radoi{\v{c}}i{\'c}, R.}
\newblock Counting 1324-avoiding permutations.
\newblock {\em Electron. J. Combin. 9}, 2 (2002/03), Research paper 13, 9 pp.
  (electronic).

\bibitem{maximillian}
{\sc Murphy, M.~M.}
\newblock {\em Restricted permutations, antichains, atomic classes, and stack
  sorting}.
\newblock PhD thesis, Univ. of St. Andrews, 2002.

\bibitem{nz:copies}
{\sc Noonan, J., and Zeilberger, D.}
\newblock The enumeration of permutations with a prescribed number of
  ``forbidden'' patterns.
\newblock {\em Adv. in Appl. Math. 17}, 4 (1996), 381--407.

\bibitem{regev}
{\sc Regev, A.}
\newblock Asymptotic values for degrees associated with strips of {Y}oung
  diagrams.
\newblock {\em Adv. in Math. 41}, 2 (1981), 115--136.

\bibitem{sw:comp}
{\sc Savage, C.~D., and Wilf, H.~S.}
\newblock Pattern avoidance in compositions and multiset permutations.
\newblock
  \href{http://front.math.ucdavis.edu/math.CO/0504310}{arXiv:math.CO/0504310}.

\bibitem{smith:sorting}
{\sc Smith, R.}
\newblock Comparing algorithms for sorting with {$t$} stacks in series.
\newblock {\em Ann. Comb. 8}, 1 (2004), 113--121.

\bibitem{stankova:len4}
{\sc Stankova, Z.}
\newblock Classification of forbidden subsequences of length {$4$}.
\newblock {\em European J. Combin. 17}, 5 (1996), 501--517.

\bibitem{sw:wilf}
{\sc Stankova, Z., and West, J.}
\newblock A new class of {W}ilf-equivalent permutations.
\newblock {\em J. Algebraic Combin. 15}, 3 (2002), 271--290.

\bibitem{stankova:fs}
{\sc Stankova, Z.~E.}
\newblock Forbidden subsequences.
\newblock {\em Discrete Math. 132}, 1-3 (1994), 291--316.

\bibitem{stanley:prec}
{\sc Stanley, R.~P.}
\newblock Differentiably finite power series.
\newblock {\em European J. Combin. 1}, 2 (1980), 175--188.

\bibitem{tarjan:sorting}
{\sc Tarjan, R.}
\newblock Sorting using networks of queues and stacks.
\newblock {\em J. Assoc. Comput. Mach. 19\/} (1972), 341--346.

\bibitem{wilfplus:prog}
{\sc Vatter, V.~R.}
\newblock \texttt{WILFPLUS}.
\newblock Available at
  \url{http://www.math.rutgers.edu/~vatter/programs/wilfplus/}.

\bibitem{west:cat}
{\sc West, J.}
\newblock Generating trees and the {C}atalan and {S}chr\"oder numbers.
\newblock {\em Discrete Math. 146}, 1-3 (1995), 247--262.

\bibitem{z:holo1}
{\sc Zeilberger, D.}
\newblock A holonomic systems approach to special functions identities.
\newblock {\em J. Comput. Appl. Math. 32}, 3 (1990), 321--368.

\end{thebibliography}
\end{small}

\end{document}